\documentclass[11pt,reqno]{amsart}
\usepackage{amsfonts}
\usepackage{amsmath,amssymb,graphicx}

\newtheorem{thm}{Theorem}[section]
\newtheorem{cor}[thm]{Corollary}
\newtheorem{lem}[thm]{Lemma}
\newtheorem{prop}[thm]{Proposition}
\theoremstyle{definition}
\newtheorem{defn}[thm]{Definition}
\theoremstyle{remark}
\newtheorem{rem}[thm]{Remark}
\numberwithin{equation}{section}

\begin{document}
\title[Doob's maximal identity and expansions of filtrations]{ Doob's
maximal identity, Multiplicative decompositions and enlargements of
filtrations }
\author{Ashkan Nikeghbali}
\address{Laboratoire de Probabilit\'es et Mod\`eles Al\'eatoires,
Universit\'e Pierre et Marie Curie, et CNRS UMR 7599, 175 rue du
Chevaleret F-75013 Paris, France.} \curraddr{American Institute of
Mathematics
360 Portage Ave Palo Alto, CA 94306-2244 \\
and University of Rochester} \email{ashkan@aimath.org}
\author{Marc Yor}
\address{Laboratoire de Probabilit\'es et Mod\`eles Al\'eatoires,
Universit\'e Pierre et Marie Curie, et CNRS UMR 7599, 175 rue du
Chevaleret F-75013 Paris, France.} 
\date{\today}
\dedicatory{In the memory of J.L. Doob}
\begin{abstract}
In the theory of progressive enlargements of filtrations, the
supermartingale $Z_{t}=\mathbf{P}\left( g>t\mid
\mathcal{F}_{t}\right) $ associated with an honest time $g$, and its
additive (Doob-Meyer) decomposition, play an essential role. In this
paper, we propose an alternative approach, using a multiplicative
representation for the supermartingale $Z_{t}$, based on Doob's
maximal identity. We thus give new examples of progressive
enlargements. Moreover, we give, in our setting, a proof of the
decomposition formula for martingales , using initial enlargement
techniques, and use it to obtain some path decompositions given the
maximum or  minimum of some processes.
\end{abstract}

\maketitle





\section{Introduction}

Let $\left( \Omega ,\mathcal{F},\left( \mathcal{F}_{t}\right) _{t\geq 0},%
\mathbf{P}\right) $ be a filtered probability space satisfying the
usual hypotheses (right continuous and complete). Given the end $L$\
of an $\left(
\mathcal{F}_{t}\right) $\ predictable set $\Gamma $, i.e\textbf{\ }%
\begin{equation*}
L=\sup \left\{ t:\left( t,\omega \right) \in \Gamma \right\} ,
\end{equation*}%
(these times are also refered to as honest times), M. Barlow
(\cite{barlow})
and Jeulin and Yor (\cite{yorjeulin}) have shown that the supermartingale:%
\begin{equation*}
Z_{t}^{L}=\mathbf{P}\left( L>t\mid \mathcal{F}_{t}\right) ,
\end{equation*}%
chosen to be c\`{a}dl\`{a}g, plays an essential role in the enlargement
formulae with respect to $L$, i.e: in expressing a general $\left( \mathcal{F%
}_{t}\right) $ martingale $\left( M_{t}\right) $ as a semimartingale in $%
\left( \mathcal{F}_{t}^{L}\right) _{t\geq 0}$, the smallest filtration which
contains $\left( \mathcal{F}_{t}\right) $, and makes $L$\ a stopping time.
This enlargement formula is:%
\begin{equation}
M_{t}=\widetilde{M}_{t}+\int_{0}^{t\wedge L}\frac{d<M,Z>_{s}}{Z_{s_{-}}}%
+\int_{L}^{t}\frac{d<M,1-Z>_{s}}{1-Z_{s_{-}}%
},  \label{grossform}
\end{equation}
where $\left( \widetilde{M}_{t}\right) _{t\geq 0}$ denotes an
$\left( \left( \mathcal{F}_{t}^{L}\right) ,\mathbf{P}\right) $ local
martingale. Hence it is important to dispose of an explicit formula
for $\left( Z_{t}^{L}\right) _{t\geq 0}$. In the literature about
progressive enlargements of filtrations, not so many examples are
fully developed (see e.g. for example \cite{zurich},
\cite{jeulinyor} or \cite{jeulin}); indeed, the computation of
$\left( Z_{t}^{L}\right) $\ is sometimes difficult. Moreover, the
examples are developed essentially in the Brownian setting, where as
we shall see, $\left( Z_{t}^{L}\right) $\ is continuous, and no
examples of discontinuous $\left( Z_{t}^{L}\right)'s$ are
given.\bigskip

In this paper, we first consider a special family of honest times
$g$, and then we later prove that this family is generic in the
sense that every honest time is in fact of this form (under some
reasonable assumptions).

More precisely, we consider the following class of local
martingales.
\begin{defn}
We say that an $\left( \mathcal{F}_{t}\right) $ local martingale
$\left( N_{t}\right) $ belongs to the class
$\left(\mathcal{C}_{0}\right)$, if it is strictly positive, with no
positive jumps, and $\lim_{t\rightarrow\infty}N_{t}=0$.
\end{defn}
\begin{rem}
Let $\left( N_{t}\right) $ be a local martingale of class
$\left(\mathcal{C}_{0}\right)$. Then: $$S_{t}\equiv \sup_{s\leq
t}N_{s},$$its supremum process, is continuous. This property is
essential in our paper. Hence, most of the results we shall state
remain valid for positive local martingales, which go to zero at
infinity, and whose suprema are continuous.
\end{rem}
We associate with a local martingale of class $\left(\mathcal{C}_{0}\right)$, the supermartingale $\left( \frac{N_{t}}{S_{t}}\right) _{t\geq 0}$%
, and the random time $g$ defined as:
\begin{eqnarray*}
g &\equiv &\sup \left\{ t\geq 0:\quad N_{t}=S_{\infty }\right\} \\
&=&\sup \left\{ t\geq 0:\quad S_{t}-N_{t}=0\right\}.
\end{eqnarray*}

In Section 2, we prove that the associated supermartingale $Z$
satisfies:
\begin{equation}
Z_{t}\equiv \mathbf{P}\left( g>t\mid \mathcal{F}_{t}\right) =\frac{N_{t}}{%
S_{t}},  \label{decomult}
\end{equation}%
and then give the decomposition formula (\ref{grossform}) in terms
of the local martingale $\left( N_{t}\right) $. This will provide us
with some new, and explicit examples of such supermartingales
$\left(Z_{t}\right)$ which are discontinuous. We also establish some
relationship between the multiplicative representation
(\ref{decomult}) and the Doob-Meyer (additive) decomposition of
$\left( Z_{t}\right) $.

In Section 3, we study the problem of the initial enlargement of $\left(\mathcal{F}%
_{t}\right)$ with the variable $S_{\infty}$, and then give a new proof of (\ref%
{grossform}).

In Section 4, we show that the formula (\ref{decomult}) is in fact very
general. More precisely, for any end of a predictable set $L$, under the
assumptions \textbf{(CA)}:

\begin{itemize}
\item all $\left( \mathcal{F}_{t}\right) $-martingales are \textbf{\underline{c}}%
ontinuous (e.g: the Brownian filtration);

\item $L$ \textbf{\underline{a}}voids every $\left( \mathcal{F}_{t}\right)$ -stopping
time $T$, i.e. $P\left[ L=T\right] =0$,
\end{itemize}
the supermartingale $Z_{t}^{L}=\mathbf{P}\left( L>t\mid \mathcal{F}%
_{t}\right) $ may be represented as (\ref{decomult}).

In Section 5, we give some new examples of enlargements of
filtrations. Moreover, as an illustration of our approach and the
method of enlargements of filtrations, we recover and complete some
known results of D. Williams (\cite{williams2}) about path
decompositions of some diffusion processes, given their minima. We
add a new fragment in these path decompositions, by introducing a
new family of random times, as defined in \cite{ANMY} and called
pseudo-stopping times, which generalize the fundamental notion of
stopping times, introduced by J.L. Doob. We take this opportunity to
quote two passages, resp. in the appendix of Meyer's book (1966):
\begin{quote}  Les temps d'arr\^{e}t ont \'{e}t\'{e}
utilis\'{e}s, sans  d\'{e}finition formelle, depuis le d\'{e}but de
la th\'{e}orie des processus. La notion appara\^{\i}t tout \`{a}
fait clairement pour la premi\`{e}re fois chez Doob en 1936.
\end{quote}and in Dellacherie-Meyer's book, volume I
(\cite{dellachmeyer}), p.184: 0194 \begin{quote} Il a sans doute
fallu autant de g\'{e}nie aux cr\'{e}ateurs du  calcul
diff\'{e}rentiel pour expliciter la notion si simple de
d\'{e}riv\'{e}e, qu'\`{a} leurs successeurs pour faire tout le
reste. L'invention des temps d'arr\^{e}t par Doob est tout \`{a}
fait comparable. \end{quote}

\section{A multiplicative representation formula}
\subsection{Doob's maximal identity}

Let $\left( N_{t}\right) _{t\geq 0}$ be a local martingale which
belongs to the class $(\mathcal{C}_{0})$, with $N_{0}=x$.
Let $S_{t}=\sup_{s\leq t}N_{s}$. We consider:%
\begin{eqnarray}
g &=&\sup \left\{ t\geq 0:\quad N_{t}=S_{\infty }\right\}  \notag \\
&=&\sup \left\{ t\geq 0:\quad S_{t}-N_{t}=0\right\} .  \label{defdeg}
\end{eqnarray}

To establish our main proposition, we shall need the following
variant of Doob's maximal inequality, which we call Doob's maximal
identity:

\begin{lem}[Doob's maximal
identity]
\label{maxeq} For any $a>0$, we have:%
\begin{enumerate}
\item
\begin{equation} \mathbf{P}\left( S_{\infty }>a\right) =\left(
\frac{x}{a}\right) \wedge 1. \label{loimax}
\end{equation}Hence, $\dfrac{x}{S_{\infty }}$ is a uniform random variable on $%
\left(0,1\right)$.
\item For any stopping time $T$:%
\begin{equation}
\mathbf{P}\left( S^{T}>a\mid \mathcal{F}_{T}\right) =\left( \frac{N_{T}}{a}%
\right) \wedge 1 ,  \label{loimaxcond}
\end{equation}%
where
\begin{equation*}
S^{T}=\sup_{u\geq T}N_{u}.
\end{equation*}%
Hence $\dfrac{N_{T}}{S^{T}}$ is also a uniform random variable on
$\left(0,1\right)$, independent of $\mathcal{F}_{T}$.
\end{enumerate}
\end{lem}

\begin{proof}
Formula (\ref{loimaxcond}) is a consequence of (\ref{loimax}) when applied
to the martingale $\left( N_{T+u}\right) _{u\geq 0}$ and the filtration $%
\left( \mathcal{F}_{T+u}\right) _{u\geq 0}$. Formula (\ref{loimax})
itself is obvious when $a\leq x$, and for $a>x$, it is obtained by
applying Doob's optional stopping theorem to the local martingale
$\left( N_{t\wedge T_{a}}\right) $, where $T_{a}=\inf \left\{ u\geq
0:\text{ }N_{u}>a\right\} $.
\end{proof}

The next proposition gives an explicit formula for $Z_{t}\equiv \mathbf{P}%
\left( g>t\mid \mathcal{F}_{t}\right) $, in terms of the local martingale $%
\left( N_{t}\right) $. Without loss of generality, \textbf{we assume
from now on that $\mathbf{x=1}$}. Indeed, if $N_{0}=x$, we consider
the local martingale $\left( \frac{N_{t}}{x}\right) $ which starts
at $1$.

\begin{prop}\label{applicationmax}
\begin{enumerate}
\item In our setting, the formula:%
\begin{equation*}
Z_{t}=\frac{N_{t}}{S_{t}},\text{ }t\geq 0
\end{equation*}%
holds.

\item The Doob-Meyer additive decomposition of $\left( Z_{t}\right) $\ is:%
\begin{equation}
Z_{t}=\mathbf{E}\left[ \log S_{\infty }\mid \mathcal{F}_{t}\right] -\log
\left( S_{t}\right) .  \label{DB}
\end{equation}
\end{enumerate}
\end{prop}

\begin{proof}
We first note that:%
\begin{eqnarray*}
\left\{ g>t\right\}  &=&\left\{ \exists \text{ }u>t:\text{ }%
S_{u}=N_{u}\right\}  \\
&=&\left\{ \exists \text{ }u>t:\text{ }S_{t}\leq N_{u}\right\}  \\
&=&\left\{ \sup_{u\geq t}N_{u}\geq S_{t}\right\} .
\end{eqnarray*}%
Hence, from (\ref{loimaxcond}), we get: $\mathbf{P}\left( g>t\mid \mathcal{F}%
_{t}\right) =\frac{N_{t}}{S_{t}}$.

To establish (\ref{DB}), we develop $\left( \frac{N_{t}}{S_{t}}\right) $\
thanks to Ito's formula, to obtain:%
\begin{equation*}
Z_{t}=1+\int_{0}^{t}\frac{1}{S_{s}}dN_{s}-\int_{0}^{t}\frac{N_{s}}{\left(
S_{s}\right) ^{2}}dS_{s}.
\end{equation*}%
Now, we remark that the measure $dS_{s}$\ is carried by the set $\left\{ s:%
\text{ }Z_{s}=1\right\} $; hence:%
\begin{eqnarray*}
Z_{t} &=&1+\int_{0}^{t}\frac{1}{S_{s}}dN_{s}-\int_{0}^{t}\frac{1}{S_{s}}%
dS_{s} \\
\dfrac{N_{t}}{S_{t}}&=&1+\int_{0}^{t}\frac{1}{S_{s}}dN_{s}-\log
\left( S_{t}\right) .
\end{eqnarray*}%
From the unicity of the Doob-Meyer decomposition, $\log \left(
S_{t}\right) $ is the predictable increasing part of $\left(
Z_{t}\right) $\ whilst $\left(
\int_{0}^{t}\frac{1}{S_{s}}dN_{s}\right) $\ is its martingale part. As $%
\left( Z_{t}\right) $\ is of class $\left( D\right) $, $\left( \int_{0}^{t}%
\frac{1}{S_{s}}dN_{s}\right) $\ is a uniformly integrable martingale. Now,
let $t\rightarrow \infty $: as $Z_{\infty }=0$, $\log S_{\infty
}=1+\int_{0}^{\infty }\frac{1}{S_{s}}dN_{s}$ and thus:
\begin{equation}\label{qqrrr}
1+\int_{0}^{t}\frac{1}{S_{s}}dN_{s}=\mathbf{E}\left[ \log S_{\infty
}\mid \mathcal{F}_{t}\right] ,
\end{equation}%
which proves (2).
\end{proof}
\begin{rem}
It is well known, and it follows from (\ref{DB}), that the
martingale in (\ref{qqrrr}) is in fact in BMO.
\end{rem}
\begin{cor}
Assuming that all $\left(\mathcal{F}_{t}\right)$ martingales are
continuous, the following hold:
\begin{enumerate}
\item $\log \left( S_{t}\right) $ is the dual predictable projection
of $\mathbf{1}_{\left\{ g\leq t\right\} }$: for any positive
predictable process $\left(k_{s}\right)$,
$$\mathbf{E}\left(k_{g}\right)=\mathbf{E}\left(\int_{0}^{\infty}k_{s}\dfrac{dS_{s}}{S_{s}}\right);$$
\item The random time $g$ is honest and avoids any $\left( \mathcal{F}_{t}\right) $%
stopping time $T$\textit{, i.e. }$P\left[ g=T\right] =0$.
\end{enumerate}
\end{cor}
\begin{proof}
Under our assumptions, the predictable and optional sigma algebras
are equal. Thus, it suffices to prove that $g$ avoids stopping
times, the other assertions being obvious. Since $\log \left(
S_{t}\right) $ is the dual predictable projection of
$\mathbf{1}_{\left\{ g\leq t\right\} }$ and is continuous, then for
any $\left( \mathcal{F}_{t}\right) $ stopping time $T$,
\begin{equation*}
\mathbf{E}\left[ \mathbf{1}_{\left\{ g=T\right\}
}\right]=\mathbf{E}\left[\left(\Delta \log \left(
S_{\bullet}\right)\right)_{T}\right]=0.
\end{equation*}%
Thus we get $P\left( g=T\right) =0$.
\end{proof}

\bigskip

We can now write the formula (\ref{grossform}) in terms of the martingale $%
\left( N_{t}\right) $.

\begin{prop}
Let $\left( X_{t}\right) _{t\geq 0}$ be a local $\left( \mathcal{F}%
_{t}\right) $ martingale. Then, $X$\ has the following decomposition as a
semimartingale in $\left( \mathcal{F}_{t}^{g}\right) $:%
\begin{equation*}
X_{t}=\widetilde{X}_{t}+\int_{0}^{t\wedge g}\frac{d<X,N>_{s}}{N_{s-}}%
-\int_{g}^{t}\frac{d<X,N>_{s}}{S _{\infty}-N_{s-}}
\end{equation*}%
where $\left( \widetilde{X}_{t}\right) $ is an $\left( \mathcal{F}%
_{t}^{g}\right) $\ local martingale.
\end{prop}

\begin{proof}
This is a consequence of formula (\ref{grossform}) and Proposition
\ref{applicationmax}.
\end{proof}

We shall now give a relationship between $\left( S_{t}\right) $\ and $%
\mathbf{E}\left[ \log S_{\infty }\mid \mathcal{F}_{t}\right] $. For
this, we shall need the following easy extension of Skorokhod's
reflection lemma (see \cite{Mckean}, p.72):
\begin{lem}\label{lemmreflection}
Let $y$ be a real-valued c\`{a}dl\`{a}g function on
$\left[0,\infty\right)$, such that $y$ has no negative jumps, and
$y(0)=0$. Then, there exists a unique pair $\left(z,a\right)$ of
functions on $\left[0,\infty\right)$ such that:
\begin{enumerate}
\item z=y+a
\item z is positive, c\`{a}dl\`{a}g and has no negative jumps,
\item a is increasing, continuous, vanishing at zero and the
corresponding measure $da_{s}$ is carried by
$\left\{s:\;z(s)=0\right\}$.
\end{enumerate} The function $a$ is moreover given by
$$a(t)=\sup_{s\leq t}\left(-y(s)\right).$$
\end{lem}

\begin{prop}
\label{sko}With
\begin{equation*}
\mu _{t}=\mathbf{E}\left[ \log S_{\infty }\mid
\mathcal{F}_{t}\right] ,
\end{equation*}%
we have:%
\begin{equation*}
\log \left( S_{t}\right) =\sup_{s\leq t}\mu _{s}-1\equiv \overline{\mu }%
_{t}-1,
\end{equation*}%
or equivalently:%
\begin{equation*}
S_{t}=\exp \left( \overline{\mu }_{t}-1\right)
\end{equation*}
\end{prop}

\begin{proof}
From (\ref{DB}), we can write:%
\begin{equation*}
1-Z_{t}=\left( 1-\mu _{t}\right) +\log \left( S_{t}\right) .
\end{equation*}%
From Lemma \ref{lemmreflection}, we deduce that
\begin{equation*}
\log \left( S_{t}\right) =\sup_{s\leq t}\mu _{s}-1.
\end{equation*}
\end{proof}
\subsection{Some hidden Az\'{e}ma-Yor martingales}
We shall now associate with the two dimensional process
\begin{equation*}
\left( \log \left( S_{t}\right) ,\text{ }Z_{t}\right) _{t\geq 0}
\end{equation*}%
a family of martingales reminiscent of Az\'{e}ma-Yor martingales
(see, e.g., \cite{AY}) which we shall now discuss. In fact, once
again, we have to introduce a slightly generalized version of what
are usually called Az\'{e}ma-Yor martingales. Indeed, these
martingales were originally defined for continuous local martingales
(see \cite{revuzyor}, Chapter VI), while we would like to define
them for local martingales without positive jumps. This extension
can be dealt with the following balayage argument:
\begin{lem}
Let $Y=M+A$ be a special semimartingale, where $M$ is a
c\`{a}dl\`{a}g local martingale, and $A$ a continuous increasing
process. Set $H=\left\{t:\;Y_{t}=0\right\}$, and define $g_{t}\equiv
\sup\left\{s<t:\;Y_{s}=0\right\}$. Then, for any locally bounded
predictable process $\left(k_{t}\right)$, $\left(k_{g_{t}}\right)$
is predictable and
\begin{equation}\label{balay}
k_{g_{t}}Y_{t}=k_{0}Y_{0}+\int_{0}^{t}k_{g_{s}}dY_{s}.
\end{equation}
\end{lem}
\begin{proof}
The proof is the same as the proof for continuous semimartingales.
The reader can refer to \cite{delmaismey}, p.144, for even more
general versions of the balayage formula.
\end{proof}Now, we can state the following generalization of the
classical Az\'{e}ma-Yor martingales:
\begin{prop}\label{azemayorgeneralisee}
Let $\left(N_{t}\right)_{t\geq 0}$ be a local martingale such that
its supremum process $\left(S_{t}\right)$ is continuous (this is the
case if $N_{t}$ is in the class $\mathcal{C}_{0}$). Let $f$ be a
locally bounded Borel function and define
$F\left(x\right)=\int_{0}^{x}dyf\left(y\right)$. Then, $X_{t}\equiv
F\left(S_{t}\right)-f\left(S_{t}\right)\left(S_{t}-N_{t}\right)$ is
a local martingale and:
\begin{equation}  \label{ayor}
F\left(S_{t}\right)-f\left(S_{t}\right)\left(S_{t}-N_{t}\right)=%
\int_{0}^{t}f\left(S_{s}\right)dN_{s}+F\left(S_{0}\right),
\end{equation}
\end{prop}
\begin{proof}
In (\ref{balay}), take $k_{t}\equiv f\left(S_{t}\right)$, and
$Y_{t}\equiv S_{t}-N_{t}$. Then, we have:
\begin{equation*}
f\left(S_{g_{t}}\right)\left(S_{t}-N_{t}\right)=\int_{0}^{t}f\left(S_{g_{s}}%
\right)d\left(S_{s}-N_{s}\right).
\end{equation*}
But $S_{g_{t}}=S_{t}$, hence:
\begin{equation*}
F\left(S_{t}\right)-f\left(S_{t}\right)\left(S_{t}-N_{t}\right)=%
\int_{0}^{t}f\left(S_{s}\right)dN_{s}+F\left(S_{0}\right).
\end{equation*}%
In conclusion, for any locally bounded function $f$,
\begin{equation*}
F\left(S_{t}\right)-f\left(S_{t}\right)\left(S_{t}-N_{t}\right)=%
\int_{0}^{t}f\left(S_{s}\right)dN_{s}+F\left(S_{0}\right),
\end{equation*}
is a local martingale. \end{proof}

\begin{rem}
Although very simple, these martingales played an essential role in
the resolution by Az\'{e}ma and Yor of Skorokhod's embedding problem
(see \cite{revuzyor}, chapter VI for more details and references).
\end{rem}
\begin{rem}
In \cite{laurentyor}, a special case of Proposition
\ref{azemayorgeneralisee}, for spectrally negative L\'{e}vy
martingales is obtained by different means.
\end{rem}
Now, we associate with the two dimensional process $\left( \log
\left( S_{t}\right) ,\text{ }Z_{t}\right) _{t\geq 0}$, a canonical
family of local martingales which are in fact of the form
(\ref{ayor}).

\begin{prop}
Let $f$ be a locally bounded and Borel function, and let $%
F\left(x\right)=\int_{0}^{x}dyf\left(y\right)$.
\begin{enumerate}
\item The following processes are local martingales:%
\begin{equation}
F\left( \log \left( S_{t}\right) \right) -f\left( \log \left( S_{t}\right)
\right) \left( 1-Z_{t}\right) ,\text{ }t\geq 0.  \label{azemayordeg}
\end{equation}

\item Denoting $K\left( x\right) =F\left( x-1\right) $ and $k\left( x\right)
=f\left( x-1\right) $, then the local martingales in (\ref{azemayordeg}) are
seen to be equal to:%
\begin{equation}
K\left( \overline{\mu }_{t}\right) -k\left( \overline{\mu
}_{t}\right) \left( \overline{\mu }_{t}-\mu _{t}\right) ,\text{
}t\geq 0.  \label{f}
\end{equation}
\end{enumerate}
\end{prop}

\begin{proof}
(1). The fact that (\ref{azemayordeg}) defines a local martingale
may be seen as an application of Ito's lemma (when $f$ is regular),
followed by a monotone class argument.

(2). Formula (\ref{f}) is obtained by a trivial change of variables,
and the fact that: $1-Z_{t}=\overline{\mu }_{t}-\mu _{t}$, which was
derived in Proposition \ref{sko}.
\end{proof}

\begin{rem}
Similar formulas are derived in \cite{ANMYII} from different
considerations.
\end{rem}
\section{Initial expansion with $S_{\infty}$ and enlargement formulae}

In this Section, we shall deal with the question of initial
enlargement of
the filtration $\left(\mathcal{F}_{t}\right)$ with the variable $S_{\infty}$%
. This problem cannot be dealt with the powerful enlargement theorem
of Jacod (see \cite{jeulinyor}), but can be treated by a careful
combination of different propositions in \cite{jeulin}. However, we
shall give a simple proof which can also be adapted to deal with
some other
situations. Eventually, we will use our result about the initial expansion of $%
\left(\mathcal{F}_{t}\right)$ with the variable $S_{\infty}$ to
recover formula (\ref{grossform}).

Let us define the new filtration
\begin{equation*}
\mathcal{F}_{t}^{\sigma\left(S_{\infty}\right)}\equiv
\bigcap_{\varepsilon>0}\left(\mathcal{F}_{t+\varepsilon}\vee
\sigma\left(S_{\infty}\right)\right),
\end{equation*}%
which satisfies the usual assumptions. The new information $%
\sigma\left(S_{\infty}\right)$ is brought in at the origin of time
and $g$ is a stopping time for this larger filtration. More
precisely:
\begin{lem}\label{lemminclusion}
The following hold:
\begin{enumerate}
\item \begin{equation*}
    g=\inf\left\{t:\;N_{t}=S_{\infty}\right\};
\end{equation*} and hence $g$ is an
$\left(\mathcal{F}_{t}^{\sigma\left(S_{\infty}\right)}\right)$
stopping time.
\item Consequently: \begin{equation*} \mathcal{F}_{t}^{g}\subset
\mathcal{F}_{t}^{\sigma\left(S_{\infty}\right)}.
\end{equation*}
\end{enumerate}
\end{lem}

\begin{proof}
$(1)$ The measure $dS_{t}$ is carried by the set
$\left\{t:\;N_{t}=S_{t}\right\}$. As
$g=\sup\left\{t:\;N_{t}=S_{t}\right\}$, the process
$\left(S_{t}\right)$ does not grow after $g$, which also satisfies:
$$g=\inf\left\{t:\;S_{t}=S_{\infty}\right\};$$hence $g$ is an
$\left(\mathcal{F}_{t}^{\sigma\left(S_{\infty}\right)}\right)$
stopping time.

$(2)$ It is obvious.
\end{proof}

Now we introduce some standard terminology.

\begin{defn}
We shall say that the pair of filtrations $\left(\mathcal{F}_{t}, \mathcal{F}%
_{t}^{\sigma\left(S_{\infty}\right)}\right)$ satisfies the $%
\left(H^{\prime}\right)$ hypothesis if every $\left(\mathcal{F}_{t}\right)$
(semi)martingale is a $\left(\mathcal{F}_{t}^{\sigma\left(S_{\infty}\right)}%
\right)$ semimartingale.
\end{defn}

We shall now show that the pair of filtrations $\left(\mathcal{F}_{t},
\mathcal{F}_{t}^{\sigma\left(S_{\infty}\right)}\right)$ satisfies the $%
\left(H^{\prime}\right)$ hypothesis and give the decomposition of a $\left(%
\mathcal{F}_{t}\right)$ local martingale in $\left(\mathcal{F}%
_{t}^{\sigma\left(S_{\infty}\right)}\right)$. For this, we need to know the
conditional law of $S_{\infty}$ given $\mathcal{F}_{t}$.

\begin{prop}
For any Borel bounded or positive function $f$, we have:
\begin{eqnarray}
\mathbf{E}\left(f\left(S_{\infty}\right)|\mathcal{F}_{t}\right) &=&
f\left(S_{t}\right)\left(1-\dfrac{N_{t}}{S_{t}}\right)+%
\int_{0}^{N_{t}/S_{t}}dxf\left(\dfrac{N_{t}}{x}\right)  \label{grosavecs} \\
&=& f\left(S_{t}\right)\left(1-\dfrac{N_{t}}{S_{t}}\right)+N_{t}%
\int_{S_{t}}^{\infty}dx\frac{f\left(x\right)}{x^{2}}.  \notag
\end{eqnarray}
\end{prop}

\begin{proof}
The proof is based on Lemma \ref{maxeq}; in the following, $U$ is a random
variable, which follows the standard uniform law and which is independent of
$\mathcal{F}_{t}$.
\begin{eqnarray*}
\mathbf{E}\left(f\left(S_{\infty}\right)|\mathcal{F}_{t}\right) &=& \mathbf{E%
}\left(f\left(S_{t}\vee S^{t}\right)|\mathcal{F}_{t}\right) \\
&=& \mathbf{E}\left(f\left(S_{t}\right)\mathbf{1}_{\left\{S_{t}\geq
S^{t}\right\}}|\mathcal{F}_{t}\right)+\mathbf{E}\left(f\left(S^{t}\right)%
\mathbf{1}_{\left\{S_{t}< S^{t}\right\}}|\mathcal{F}_{t}\right) \\
&=& f\left(S_{t}\right)\mathbf{P}\left(S_{t}\geq S^{t}|\mathcal{F}%
_{t}\right)+ \mathbf{E}\left(f\left(S^{t}\right)\mathbf{1}_{\left\{S_{t}<
S^{t}\right\}}|\mathcal{F}_{t}\right) \\
&=& f\left(S_{t}\right)\mathbf{P}\left(U\leq \dfrac{N_{t}}{S_{t}}|\mathcal{F}_{t}\right)+%
\mathbf{E}\left(f\left(\dfrac{N_{t}}{U}\right)\mathbf{1}_{\left\{U<\frac{%
N_{t}}{S_{t}}\right\}}|\mathcal{F}_{t}\right) \\
&=& f\left(S_{t}\right)\left(1-\dfrac{N_{t}}{S_{t}}\right)+%
\int_{0}^{N_{t}/S_{t}}dxf\left(\dfrac{N_{t}}{x}\right).
\end{eqnarray*}%
A straightforward change of variable in the last integral also gives:
\begin{equation*}
\mathbf{E}\left(f\left(S_{\infty}\right)|\mathcal{F}_{t}\right)=f\left(S_{t}%
\right)\left(1-\dfrac{N_{t}}{S_{t}}\right)+N_{t}\int_{S_{t}}^{\infty}dy\frac{%
f\left(y\right)}{y^{2}}.
\end{equation*}
\end{proof}

One may now ask if $\mathbf{E}\left(f\left(S_{\infty}\right)|\mathcal{F}%
_{t}\right)$ is of the form (\ref{ayor}). The answer to this question is
positive. Indeed:
\begin{eqnarray*}
\mathbf{E}\left(f\left(S_{\infty}\right)|\mathcal{F}_{t}\right)&=&
f\left(S_{t}\right)\left(1-\dfrac{N_{t}}{S_{t}}\right)+N_{t}\int_{S_{t}}^{%
\infty}dy\frac{f\left(y\right)}{y^{2}} \\
&=& S_{t}\int_{S_{t}}^{\infty}dy\frac{f\left(y\right)}{y^{2}}%
-\left(S_{t}-N_{t}\right)\left(\int_{S_{t}}^{\infty}dy\frac{f\left(y\right)}{%
y^{2}}-\dfrac{f\left(S_{t}\right)}{S_{t}}\right).
\end{eqnarray*}%
Hence,
\begin{equation*}
\mathbf{E}\left(f\left(S_{\infty}\right)|\mathcal{F}_{t}\right)=H\left(1%
\right)+H\left(S_{t}\right)-h\left(S_{t}\right)\left(S_{t}-N_{t}\right),
\end{equation*}
with
\begin{equation*}
H\left(x\right)=x\int_{x}^{\infty}dy\frac{f\left(y\right)}{y^{2}},
\end{equation*}
and
\begin{equation*}
h\left(x\right)=h_{f}\left(x\right)\equiv\int_{x}^{\infty}dy\frac{f\left(y\right)}{y^{2}}-\dfrac{%
f\left(x\right)}{x}=\int_{x}^{\infty}\frac{dy}{y^{2}}\left(f\left(y\right)-f\left(x\right)\right).
\end{equation*}
Moreover, again from formula (\ref{ayor}), we have the following
representation of $\mathbf{E}\left(f\left(S_{\infty}\right)|\mathcal{F}%
_{t}\right)$ as a stochastic integral:
\begin{equation}  \label{represstoc}
\mathbf{E}\left(f\left(S_{\infty}\right)|\mathcal{F}_{t}\right)=\mathbf{E}%
\left(f\left(S_{\infty}\right)\right)+\int_{0}^{t}h\left(S_{s}\right)dN_{s}.
\end{equation}%
Let us sum up these results, introducing some notations:
\begin{eqnarray}
\lambda_{t}\left(f\right) &\equiv& \mathbf{E}\left(f\left(S_{\infty}\right)|%
\mathcal{F}_{t}\right) \\
&=& f\left(S_{t}\right)\left(1-\dfrac{N_{t}}{S_{t}}\right)+N_{t}%
\int_{S_{t}}^{\infty}dx\frac{f\left(x\right)}{x^{2}};
\end{eqnarray}%
and
\begin{equation}
\lambda_{t}\left(f\right)=\mathbf{E}\left(f\left(S_{\infty}\right)\right)+%
\int_{0}^{t}\dot{\lambda}_{s}\left(f\right)dN_{s},
\end{equation}%
where:
\begin{equation}
\dot{\lambda}_{s}\left(f\right)=h_{f}\left(S_{s}\right).
\end{equation}%
Moreover, there exist two families of random measures $%
\left(\lambda_{t}\left(dx\right)\right)_{t\geq 0}$ and $\left(\dot{\lambda}%
_{t}\left(dx\right)\right)_{t\geq 0}$, with
\begin{eqnarray}
\lambda_{t}\left(dx\right) &=& \left(1-\dfrac{N_{t}}{S_{t}}%
\right)\delta_{S_{t}}\left(dx\right)+N_{t}\mathbf{1}_{\left\{x>S_{t}\right\}}%
\dfrac{dx}{x^{2}} \\
\dot{\lambda}_{t}\left(dx\right) &=& -\dfrac{1}{S_{t}}\delta_{S_{t}}\left(dx%
\right)+\mathbf{1}_{\left\{x>S_{t}\right\}}\dfrac{dx}{x^{2}},
\end{eqnarray}
such that
\begin{eqnarray}
\lambda_{t}\left(f\right) &=& \int\lambda_{t}\left(dx\right)f\left(x\right)
\\
\dot{\lambda}_{t}\left(f\right) &=& \int\dot{\lambda}_{t}\left(dx\right)f%
\left(x\right).
\end{eqnarray}%
Eventually, we notice that there is an absolute continuity relationship
between $\lambda_{t}\left(dx\right)$ and $\dot{\lambda}_{t}\left(dx\right)$;
more precisely,
\begin{equation}
\dot{\lambda}_{t}\left(dx\right)=\lambda_{t}\left(dx\right)\rho\left(x,t%
\right),
\end{equation}%
with
\begin{equation}  \label{absolucontrel}
\rho\left(x,t\right)=\dfrac{-1}{S_{t}-N_{t}}\mathbf{1}_{\left\{S_{t}=x\right%
\}}+\dfrac{1}{N_{t}}\mathbf{1}_{\left\{S_{t}<x\right\}}.
\end{equation}%
Now, we can state the main theorem of this section.

\begin{thm}
\label{decoinitial} Let $\left(N_{t}\right)_{t\geq 0}$ be a local
martingale in the class $\mathcal{C}_{0}$ (recall $N_{0}=1$).
Then, the pair of filtrations $\left(\mathcal{F}_{t}, \mathcal{F}%
_{t}^{\sigma\left(S_{\infty}\right)}\right)$ satisfies the $%
\left(H^{\prime}\right)$ hypothesis and every $\left(\mathcal{F}_{t}\right)$
local martingale $\left(X_{t}\right)$ is an $\left(\mathcal{F}_{t}^{\sigma\left(S_{\infty}%
\right)}\right)$ semimartingale with canonical decomposition:
\begin{equation*}
X_{t}=\widetilde{X}_{t}+\int_{0}^{t}\mathbf{1}_{\left\{ g>s\right\} }\frac{%
d<X,N>_{s}}{N_{s-}}-\int_{0}^{t}\mathbf{1}_{\left\{ g\leq s\right\} }\frac{%
d<X,N>_{s}}{S_{\infty}-N_{s-}},
\end{equation*}%
where $\left( \widetilde{X}_{t}\right) $ is a $\left(\mathcal{F}%
_{t}^{\sigma\left(S_{\infty}\right)}\right)$\ local martingale.
\end{thm}
\begin{rem}
The following proof is tailored on the arguments found in
\cite{zurich}, although our framework is more general: we do not
assume that our filtration has the predictable representation
property with respect to some martingale nor that all martingales
are continuous.
\end{rem}

\begin{proof}
We can first assume that $X$ is in $\mathcal{H}^{1}$; the general
case follows by localization. Let $\Lambda_{s}$ be an
$\mathcal{F}_{s}$ measurable set, and take $t>s$. Then, for any
bounded test function $f$, $\lambda_{t}\left(f\right)$ is a bounded
martingale, hence in $BMO$, and we
have:%
\begin{eqnarray*}
\mathbf{E}\left(\mathbf{1}_{\Lambda_{s}}f\left( A_{\infty
}\right)\left(X_{t}-X_{s}\right)\right) &=& \mathbf{E}\left(\mathbf{1}%
_{\Lambda_{s}}\left(\lambda_{t}\left(f\right)X_{t}-\lambda_{s}\left(f%
\right)X_{s}\right)\right) \\
&=& \mathbf{E}\left(\mathbf{1}_{\Lambda_{s}}\left(<\lambda\left(f%
\right),X>_{t}-<\lambda\left(f\right),X>_{s}\right)\right) \\
&=& \mathbf{E}\left(\mathbf{1}_{\Lambda_{s}}\left(\int_{s}^{t}\dot{\lambda}%
_{u}\left(f\right)d<X,N>_{u}\right)\right) \\
&=& \mathbf{E}\left(\mathbf{1}_{\Lambda_{s}}\left(\int_{s}^{t}\int%
\lambda_{u}\left(dx\right)\rho\left(x,u\right)f\left(x\right)d<X,N>_{u}%
\right)\right) \\
&=& \mathbf{E}\left(\mathbf{1}_{\Lambda_{s}}\left(\int_{s}^{t}d<X,N>_{u}\rho%
\left(S_{\infty },u\right)\right)\right).
\end{eqnarray*}%
But from (\ref{absolucontrel}), we have:%
\begin{equation*}
\rho\left(S_{\infty },t\right)=\dfrac{-1}{S_{t}-N_{t}}\mathbf{1}%
_{\left\{S_{t}=S_{\infty}\right\}}+\dfrac{1}{N_{t}}\mathbf{1}_{\left\{S_{t}<S_{\infty}\right%
\}}.
\end{equation*}
It now  suffices to note (from Lemma \ref{lemminclusion}) that
$\left(S_{t}\right)$ is constant after  $g$ and $g$ is the first
time when $S_{\infty}=S_{t}$, or in other words:
\begin{equation*}
\mathbf{1}_{\left\{S_{\infty}>S_{t}\right\}}=\mathbf{1}_{\left\{g>t\right\}},%
\text{ and }\mathbf{1}_{\left\{S_{\infty}=S_{t}\right\}}=\mathbf{1}%
_{\left\{g\leq t\right\}}.
\end{equation*}%
This completes the proof.
\end{proof}

Theorem \ref%
{decoinitial} yields a new proof of the decomposition formula in the
progressive enlargement case. More precisely, we have:

\begin{cor}
\label{hyphprimepourN} The pair of filtrations $\left(\mathcal{F}_{t},%
\mathcal{F}_{t}^{g}\right)$ satisfies the $\left(H^{\prime}\right)$
hypothesis. Moreover, every $\left(\mathcal{F}_{t}\right)$ local martingale $%
X$ decomposes as:
\begin{equation*}
X_{t}=\widetilde{X}_{t}+\int_{0}^{t}\mathbf{1}_{\left\{ g>s\right\} }\frac{%
d<X,N>_{s}}{N_{s}}-\int_{0}^{t}\mathbf{1}_{\left\{ g\leq s\right\} }\frac{%
d<X,N>_{s}}{S_{\infty}-N_{s}},
\end{equation*}%
where $\left( \widetilde{X}_{t}\right) $ is a $\left(\mathcal{F}%
_{t}^{g}\right)$\ local martingale.
\end{cor}

\begin{proof}
Let $X$ be an $\left(\mathcal{F}_{t}\right)$ martingale which is in
$\mathcal{H}^{1}$; the general case follows by localization. From
Theorem \ref{decoinitial}
\begin{equation*}
X_{t}=\widetilde{X}_{t}+\int_{0}^{t}\mathbf{1}_{\left\{ g>s\right\} }\frac{%
d<X,N>_{s}}{N_{s}}-\int_{0}^{t}\mathbf{1}_{\left\{ g\leq s\right\} }\frac{%
d<X,N>_{s}}{S_{\infty}-N_{s}},
\end{equation*}%
where $\left(\widetilde{X}_{t}\right) _{t\geq 0}$ denotes an  $\left(%
\mathcal{F}_{t}^{\sigma\left(S_{\infty}\right)}\right)$ martingale.
Thus, $\left(\widetilde{X}_{t}\right) $, which is equal to:
\begin{equation*}
X_{t}-\left(\int_{0}^{t}\mathbf{1}_{\left\{ g>s\right\} }\frac{d<X,N>_{s}}{%
N_{s}}-\int_{0}^{t}\mathbf{1}_{\left\{ g\leq s\right\} }\frac{d<X,N>_{s}}{%
S_{\infty}-N_{s}},\right),
\end{equation*}
is $\left(\mathcal{F}_{t}^{g}\right)$ adapted (recall that $\mathcal{F}%
_{t}^{g}\subset \mathcal{F}_{t}^{\sigma\left(S_{\infty}\right)}$),
and hence it is an $\left(\mathcal{F}_{t}^{g}\right)$ martingale.
\end{proof}

\section{A multiplicative characterization of $Z_{t}$}

Usually, in the literature about progressive enlargements of
filtrations, it is assumed that the conditions \textbf{(CA)} are
satisfied. Now, we shall prove that under this assumption the supermartingale $%
Z_{t}^{L}=\mathbf{P}\left( L>t\mid \mathcal{F}_{t}\right) $, associated with
an honest time, can be represented as $\left( \dfrac{N_{t}}{S_{t}}\right)
_{t\geq 0}$, where $N_{t}$ is a positive local martingale. More precisely,
we have the following:

\begin{thm}
\label{multiplicatcarac} Let $L$\ be an honest time. Then, under the
conditions \textbf{(CA)}, there exists a continuous and nonnegative
local martingale $\left( N_{t}\right) _{t\geq 0}$,
with $N_{0}=1$ and $\lim_{t\rightarrow \infty }N_{t}=0$, such that:%
\begin{equation*}
Z_{t}=\mathbf{P}\left( L>t\mid \mathcal{F}_{t}\right) =\dfrac{N_{t}}{S_{t}}
\end{equation*}
\end{thm}

\begin{proof}
Under the conditions \textbf{(CA)}, $\left( Z_{t}\right) _{t\geq 0}$ is
continuous and can be written as (see \cite{azema} or \cite{delmaismey} for
details):%
\begin{equation*}
Z_{t}=M_{t}-A_{t},
\end{equation*}%
where $\left( M_{t}\right) $ and $\left( A_{t}\right) $ are continuous, $%
Z_{0}=1$ and $dA_{t}$ is carried by $\left\{ t:\text{ }Z_{t}=1\right\} $.
Then, for $t<T_{0}\equiv \inf\left\{t:\;Z_{t}=0\right\}$, we have:%
\begin{equation*}
\log \left( Z_{t}\right) =\int_{0}^{t}\frac{dM_{s}}{Z_{s}}-\frac{1}{2}%
\int_{0}^{t}\frac{d<M>_{s}}{Z_{s} ^{2}}-A_{t},
\end{equation*}%
hence:%
\begin{equation}
-\log \left( Z_{t}\right) =-\left( \int_{0}^{t}\frac{dM_{s}}{Z_{s}}-\frac{1}{%
2}\int_{0}^{t}\frac{d<M>_{s}}{Z_{s} ^{2}}\right) +A_{t}; \label{a}
\end{equation}%
and, from Skorokhod's reflection lemma, we have:
\begin{equation} \label{logito}
A_{t}=\sup_{u\leq t}\left( \int_{0}^{u}\frac{dM_{s}}{Z_{s}}-\frac{1}{2}%
\int_{0}^{u}\frac{d<M>_{s}}{Z_{s} ^{2}}\right) .
\end{equation}%
Now, combining (\ref{a}) and (\ref{logito}), we obtain:
\begin{equation*}
Z_{t}=\frac{N_{t}}{S_{t}},
\end{equation*}%
where
\begin{equation*}
N_{t}=\exp \left( \int_{0}^{t}\frac{dM_{s}}{Z_{s}}-\frac{1}{2}\int_{0}^{t}%
\frac{d<M>_{s}}{Z_{s} ^{2}}\right)
\end{equation*}%
is a local martingale starting from $1$, and
\begin{eqnarray*}
S_{t} &=&\sup_{u\leq t}\left(\exp \left( \int_{0}^{u}\frac{dM_{s}}{Z_{s}}-\frac{1}{%
2}\int_{0}^{u}\frac{d<M>_{s}}{Z_{s} ^{2}}\right)\right) \\
&=&\exp \left( \sup_{u\leq t}\left( \int_{0}^{u}\frac{dM_{s}}{Z_{s}}-\frac{1%
}{2}\int_{0}^{u}\frac{d<M>_{s}}{Z_{s} ^{2}}\right) \right) \\
&=&\exp \left( A_{t}\right) .
\end{eqnarray*}We finally note that, since $Z_{T_{0}}=0$, $\lim_{t\uparrow
T_{0}}N_{t}=0$, which allows to define $N_{t}$ for all $t\geq 0$.
\end{proof}
\begin{cor}
The supermartingale $Z_{t}=\mathbf{P}\left( L>t\mid
\mathcal{F}_{t}\right)$ admits the following additive and
multiplicative representations:
\begin{eqnarray*}
  Z_{t} &=& \dfrac{N_{t}}{S_{t}} \\
  Z_{t} &=& M_{t}-A_{t}.
\end{eqnarray*} Moreover, these two representations are related as follows:
\begin{eqnarray*}
  N_{t} &=& \exp \left( \int_{0}^{t}\frac{dM_{s}}{Z_{s}}-\frac{1}{2}\int_{0}^{t}%
\frac{d<M>_{s}}{Z_{s} ^{2}}\right) \\
  S_{t} &=& \exp\left(A_{t}\right);
\end{eqnarray*}and
\begin{eqnarray*}
  M_{t} &=& 1+\int_{0}^{t}\dfrac{dN_{s}}{S_{s}}=\mathbf{E}\left(\log S_{\infty}\mid \mathcal{F}_{t}\right), \\
  A_{t} &=& \log S_{t}.
\end{eqnarray*}
\end{cor}
\begin{proof}
It is a consequence of Proposition \ref{applicationmax} and Theorem
\ref{multiplicatcarac}.
\end{proof}
\bigskip

Now, as a consequence of Theorem \ref{multiplicatcarac}, we can recover the
enlargement formulae and the fact that the pair of filtrations $\left(%
\mathcal{F}_{t}, \mathcal{F}_{t}^{L}\right)$ satisfies the $(H^{\prime})$
hypothesis:

\begin{cor}
Let $L$\ be an honest time. Then under the conditions \textbf{(CA)},
the pair of filtrations $\left(\mathcal{F}_{t},
\mathcal{F}_{t}^{L}\right)$ satisfies the $(H^{\prime})$ hypothesis
and
every $\left(\mathcal{F}_{t}\right)$ local martingale $X$ is an $\left(%
\mathcal{F}_{t}^{L}\right)$ semimartingale with canonical decomposition:
\begin{equation*}
X_{t}=\widetilde{X}_{t}+\int_{0}^{t\wedge L}\dfrac{d<X,Z>_{s}}{Z_{s}}%
+\int_{L}^{t}\dfrac{d<X,1-Z>_{s}}{1-Z_{s}%
},
\end{equation*}
where $\left( \widetilde{X}_{t}\right) _{t\geq 0}$ denotes an
$\left( \left( \mathcal{F}_{t}^{L}\right)\right) $ local martingale.
\end{cor}

\begin{proof}
It is a combination of Theorem \ref{multiplicatcarac} and Corollary \ref%
{hyphprimepourN}.
\end{proof}
\begin{rem}
We then see that under the assumptions \textbf{(CA)}, the initial
enlargement of filtrations with $A_{\infty}$ amounts to enlarging
initially the filtration with $S_{\infty}$, the terminal value of
the supremum process of a continuous local martingale in
$\mathcal{C}_{0}$.
\end{rem}

We shall now outline another nontrivial consequence of Theorem
\ref{multiplicatcarac} here. In \cite{azemjeulknightyor}, the
authors are interested in giving explicit examples of dual
predictable projections of processes of the form $\mathbf{1}_{g\leq
t}$, where $g$ is an honest time. Indeed, these dual projections are
natural examples of increasing injective processes (see
\cite{azemjeulknightyor} for more details and references). With
Theorem \ref{multiplicatcarac}, we have a complete characterization
of such projections:
\begin{cor}
Assume the assumption \textbf{(C)} holds, and let
$\left(C_{t}\right)$ be an increasing process. Then $C$ is the dual
predictable projection of $\mathbf{1}_{g\leq t}$, for some honest
time $g$ that avoids stopping times, if and only if there exists a
continuous local martingale $N_{t}$ in the class $\mathcal{C}_{0}$
such that
$$C_{t}=\log S_{t}.$$
\end{cor}

\bigskip
The previous results can be naturally extended to the case where the
supermartingale $Z_{t}$ has only negative jumps; we gave a special
treatment under the hypothesis \textbf{(CA)} because of its
practical importance. We just give here the extension of Theorem
\ref{multiplicatcarac}; the corollaries are easily deduced.
\begin{prop}
Let $L$\ be an honest time that avoids stopping times. Assume that
$Z_{t}^{L}$ has no positive jumps. Then, there exists a local
martingale $\left( N_{t}\right) _{t\geq 0}$, in the class
$\mathcal{C}_{0}$,
with $N_{0}=1$, such that:%
\begin{equation*}
\left(Z_{t}^{L}=\right)Z_{t}=\mathbf{P}\left( L>t\mid
\mathcal{F}_{t}\right) =\dfrac{N_{t}}{S_{t}}
\end{equation*}
\end{prop}
\begin{proof}
We use the same notations as in the proof of Theorem
\ref{multiplicatcarac}. For $t<T_{0}\equiv \inf\left\{t:\;Z_{t}=0\right\}$, we have:%
\begin{equation*}
-\log \left( Z_{t}\right) =-\left(\int_{0}^{t}\left(\frac{dM_{s}}{Z_{s-}}-\frac{1}{2}%
\frac{d<M^{c}>_{s}}{Z_{s-} ^{2}}\right)+\sum_{0<s\leq t}\left(\log
\left(1+\dfrac{\Delta Z_{s}}{Z_{s-}}\right)-\dfrac{\Delta
Z_{s}}{Z_{s-}}\right)\right)+A_{t}.
\end{equation*}Now, from Lemma \ref{lemmreflection},
$$A_{t}=\sup_{s\leq t}\left(\int_{0}^{t}\left(\frac{dM_{s}}{Z_{s-}}-\frac{1}{2}%
\frac{d<M^{c}>_{s}}{Z_{s-} ^{2}}\right)+\sum_{0<s\leq t}\left(\log
\left(1+\dfrac{\Delta Z_{s}}{Z_{s-}}\right)-\dfrac{\Delta
Z_{s}}{Z_{s-}}\right)\right).$$ Now, combining the last two
equalities, we obtain:
$$Z_{t}=\dfrac{N_{t}}{S_{t}},$$where $$N_{t}=\exp\left(\int_{0}^{t}\left(\frac{dM_{s}}{Z_{s-}}-\frac{1}{2}%
\frac{d<M^{c}>_{s}}{Z_{s-} ^{2}}\right)\right)\prod_{0<s\leq
t}\left(1+\dfrac{\Delta
Z_{s}}{Z_{s-}}\right)\exp\left(-\dfrac{\Delta
Z_{s}}{Z_{s-}}\right).$$
\end{proof}

\section{Examples and applications}
In this section, we look at some specific local martingales $N_{t}$,
and use the initial enlargement formula with $S_{\infty}$, to get
some path decompositions, given the maximum or the minimum of some
stochastic processes. Our aim here is to illustrate how techniques
from enlargement of filtrations can be applied. To have a complete
description for the path decompositions, we associate with $g$ a
random time, called pseudo-stopping time, which occurs before $g$.
Eventually, we give some explicit examples of supermartingales
$Z_{t}$ with jumps.
\subsection{Pseudo-stopping times}\label{secpta}
In \cite{ANMY}, we have proposed the following generalization of
stopping times:
\begin{defn}
Let $\rho:\;(\Omega,\mathcal{F})\rightarrow\mathbf{R}_{+}$ be a
random time; $\rho$ is called a pseudo-stopping time if for every
bounded $\left(\mathcal{F}_{t}\right)$ martingale we have:
$$\mathbf{E}\left(M_{\rho}\right)=\mathbf{E}\left(M_{0}\right).$$
\end{defn}David Williams (\cite{williams}) gave the first example of such a random
time and the following systematic construction is established in
\cite{ANMY}:
\begin{prop}\label{ptaconstruction}
Let $L$ be an honest time. Then, under the conditions \textbf{(CA)},
$$\rho\equiv \sup\left\{t<L:\;Z_{t}^{L}=\inf_{u\leq
L}Z_{u}^{L}\right\},$$is a pseudo-stopping time, with
$$Z_{t}^{\rho}\equiv
\mathbf{P}\left(\rho>t\mid\mathcal{F}_{t}\right)=\inf_{u\leq
t}Z_{u}^{L},$$and $Z_{\rho}^{\rho}$ follows the uniform distribution
on $(0,1)$.
\end{prop}The following property, also proved in \cite{ANMY}, is
essential in studying path decompositions:
\begin{prop}\label{regenrative}
Let $\rho$ be a pseudo-stopping time and let $M_{t}$ be an
$\left(\mathcal{F}_{t}\right)$ local martingale. Then
$\left(M_{t\wedge \rho}\right)$ is an
$\left(\mathcal{F}_{t}^{\rho}\right)$ local martingale.
\end{prop}In our setting, Proposition \ref{ptaconstruction} gives:
\begin{prop}\label{ptamult}
Define the nonincreasing process $\left(r_{t}\right)$ by:
$$r_{t}\equiv \inf_{u\leq t}\dfrac{N_{u}}{S_{u}}.$$Then,
$$\rho\equiv \sup\left\{t<g:\;\dfrac{N_{t}}{S_{t}}=\inf_{u\leq
g}\dfrac{N_{u}}{S_{u}}\right\},$$is a pseudo-stopping time and
$r_{\rho}$ follows the uniform distribution on $(0,1)$.
\end{prop}
\subsection{Path decompositions given the maxima or the minima of a diffusion}
Now, we shall apply the techniques of enlargements of filtrations to
establish some path decompositions results. Some of the following
results have been proved by David Williams in \cite{williams2},
using different methods. Jeulin has also given a proof based on
enlargements techniques in the case of transient diffusions (see
\cite{jeulin}). Here, we complete the results of David Williams by
introducing the pseudo-stopping times $\rho$ defined in Proposition
\ref{ptamult}, and we detail some interesting examples.
\subsubsection{The killed Brownian Motion}
Let $$N_{t}\equiv B_{t},$$where $\left(B_{t}\right)_{t\geq 0}$ is a
Brownian Motion starting at $1$, and stopped at
$T_{0}=\inf\left\{t:\;B_{t}=0\right\}$. Let
$$S_{t}\equiv \sup_{s\leq t}B_{s}.$$ Let $$g=\sup\left\{t:B_{t}=S_{t}\right\}$$ and $$\rho=\sup\left\{t<g:\;\dfrac{B_{t}}{S_{t}}=\inf_{u\leq
g}\dfrac{B_{u}}{S_{u}}\right\}.$$ From Doob's maximal identity,
$S_{T_{0}}=S_{g}$ is distributed as the reciprocal of a uniform
distribution $\left(0,1\right)$, i.e. it has the density:
$\mathbf{1}_{\left[1,\infty\right)}\left(x\right)\dfrac{1}{x^{2}}$.
\begin{prop}\label{madeco}
Let $\left(B_{t}\right)_{t\geq 0}$ be a Brownian Motion starting at
$1$ and stopped when it first hits $0$. Then:
\begin{itemize}
\item $\dfrac{B_{\rho}}{S_{\rho}}$ follows the uniform law on $(0,1)$, and conditionally on $\dfrac{B_{\rho}}{S_{\rho}}=r$,
$\left(B_{t}\right)$ is a Brownian Motion up to the first time when
$B_{t}=rS_{t}$.
 \item $\left(B_{t}\right)$
is an $\left(\mathcal{F}_{t}^{g}\right)$ and
$\left(\mathcal{F}_{t}^{\sigma\left(S_{T_{0}}\right)}\right)$
semimartingale with canonical decomposition:
\begin{equation}\label{decobessel}
    B_{t}=\widetilde{B}_{t}+\int_{0}^{t\wedge g}\dfrac{ds}{B_{s}}-\int_{g}^{t\wedge
    T_{0}}\dfrac{ds}{S_{T_{0}}-B_{s}},
\end{equation}where $\left(\widetilde{B}_{t}\right)$ is an
$\mathcal{F}_{t}^{\sigma\left(S_{T_{0}}\right)}$ Brownian Motion,
stopped at $T_{0}$ and independent of $S_{T_{0}}$. Consequently, we
have the following path decomposition: conditionally on
$S_{T_{0}}=m$:
\begin{enumerate}
\item the process $\left(B_{t};\;t\leq g\right)$ is a Bessel process of
dimension $3$, started from $1$, considered up to $T_{m}$, the first
time when it hits $m$;
\item the process $\left(S_{g}-B_{g+t};\;t\leq
T_{0}-g\right)$ is a $\left(\mathcal{F}_{g+t}\right)$ three
dimensional Bessel process, started from $0$, considered up to
$T_{m}$, the first time when it hits $m$, and is independent of
$\left(B_{t};\;t\leq g\right)$.
\end{enumerate}
\end{itemize}
\end{prop}
\begin{proof}
The results concerning the decomposition until $\rho$ are
consequences of the results of Subsection \ref{secpta}. The
decomposition formula is a consequence of Theorem \ref{decoinitial}.
Since $\left(\widetilde{B}_{t}\right)$ is an
$\mathcal{F}_{t}^{\sigma\left(S_{T_{0}}\right)}$ local martingale,
with $t\wedge T_{0}$ as its bracket, it follows from L\'{e}vy's
theorem that it is an
$\mathcal{F}_{t}^{\sigma\left(S_{T_{0}}\right)}$ Brownian Motion.
Moreover, it is independent of
$\mathcal{F}_{0}^{\sigma\left(S_{T_{0}}\right)}=\sigma\left(S_{T_{0}}\right)$.
Now, conditionally on $S_{T_{0}}=m$, with
$T_{m}=\inf\left\{t:\;B_{t}=m\right\}$, $\left(B_{t}\right)$
satisfies the following stochastic differential equation:
$$B_{t}=\widetilde{B}_{t}+\int_{0}^{t\wedge
T_{m}}\frac{ds}{B_{s}}.$$Hence it is a three dimensional Bessel
process up to $T_{m}$.

It also follows from the decomposition formula that:
$$B_{g+t}=\widetilde{B}_{g+t}+\int_{0}^{g}\dfrac{ds}{B_{s}}-\int_{0}^{t\wedge
    (T_{0}-g)}\dfrac{ds}{S_{g}-B_{g+s}}.$$ This equation can also be written
    as:$$S_{g}-B_{g+t}=-\left(\widetilde{B}_{g+t}-\widetilde{B}_{g}\right)+\int_{0}^{t\wedge
    (T_{0}-g)}\dfrac{ds}{S_{g}-B_{g+s}}.$$ Now,
    $\left(\widetilde{B}_{g+t}-\widetilde{B}_{g}\right)$ is an
    $\left(\mathcal{F}_{g+t}\right)$ Brownian Motion, starting from $0$, and is independent of $\mathcal{F}_{g}$. Taking
    $\widetilde{\beta}_{t}\equiv
    -\left(\widetilde{B}_{g+t}-\widetilde{B}_{g}\right)$, which is
    also an
    $\left(\mathcal{F}_{g+t}\right)$ Brownian Motion, starting from $0$, independent of
    $\mathcal{F}_{g}$, the process $\xi_{t}\equiv S_{g}-B_{t}$
    satisfies the stochastic differential equation:
    $$\xi_{t}=\widetilde{\beta}_{t}+\int_{0}^{t\wedge
(T_{0}-g)}\frac{ds}{\xi_{s}};$$hence it is a three dimensional
Bessel process, started at $0$, and considered up to $T_{m}$, and
conditionally on $S_{g}$, is independent of $\mathcal{F}_{g}$.
\end{proof}
\subsubsection{Some recurrent diffusions}
The previous example can be generalized to a wider class of
recurrent diffusions $\left(X_{t}\right)$, satisfying the stochastic
differential equation:
\begin{equation}\label{equationrecurrence}
X_{t}=x+B_{t}+\int_{0}^{t}b\left(X_{s}\right)ds,\;x>0
\end{equation}where $\left(B_{t}\right)$ is the standard Brownian
Motion, and $b$ is a Borel integrable function which allows
existence and uniqueness for equation (\ref{equationrecurrence})
(for example $b$ bounded or Lipschitz continuous). The infinitesimal
generator $L$ of this diffusion is:
$$L=\frac{1}{2}\dfrac{d^{2}}{dx^{2}}+b\left(x\right)\dfrac{d}{dx}.$$Let $T_{0}\equiv
\inf\left\{t:\;X_{t}=0\right)$, and denote by $s$ the scale function
of $X$, which is strictly increasing and which vanishes at zero,
i.e:
$$s\left(z\right)=\int_{0}^{z}\exp\left(-2\widehat{b}\left(y\right)\right)dy,$$where$$\widehat{b}\left(y\right)=\int_{0}^{y}b\left(u\right)du.$$ Hence,
$$N_{t}\equiv\dfrac{s\left(X_{t\wedge T_{0}}\right)}{s\left(x\right)}$$ is a continuous local
martingale belonging to the class $\mathcal{C}_{0}$. If $S_{t}$
denotes the supremum process of $N_{t}$ and $\overline{X}_{t}$ the
supremum process of $X_{t}$, we
have:$$S_{t}=\dfrac{s\left(\overline{X}_{t\wedge
T_{0}}\right)}{s\left(x\right)}.$$Now,
let$$g=\sup\left\{t<T_{0}:\;X_{t}=\overline{X}_{t}\right\},$$and$$\rho=\sup\left\{t<g:\;\dfrac{X_{t}}{\overline{X}_{t}}=\inf_{u\leq
g}\dfrac{X_{u}}{\overline{X}_{u}}\right\}.$$
\begin{prop}
Let $\left(X_{t}\right)$ be a diffusion process satisfying equation
(\ref{equationrecurrence}). Then:
\begin{itemize}
\item $\dfrac{X_{\rho}}{\overline{X}_{\rho}}$ follows the uniform law on $(0,1)$, and conditionally on $\dfrac{X_{\rho}}{\overline{X}_{\rho}}=r$,
$\left(X_{t},\;t\leq\rho\right)$ is a diffusion process, up to the
first time when $X_{t}=r\overline{X}_{t}$, with the same
infinitesimal generator as $X$.
 \item $\left(X_{t}\right)$
is an $\left(\mathcal{F}_{t}^{g}\right)$ and an
$\left(\mathcal{F}_{t}^{\sigma\left(\overline{X}_{T_{0}}\right)}\right)$
semimartingale with canonical decomposition:
\begin{equation}\label{decobessel2}
    X_{t}=\widetilde{B}_{t}+\int_{0}^{t}b\left(X_{u}\right)du+\int_{0}^{t\wedge g}\dfrac{s'\left(X_{u}\right)}{s\left(X_{u}\right)}du-\int_{g}^{t\wedge
    T_{0}}\dfrac{s'\left(X_{u}\right)}{s\left(\overline{X}_{T_{0}}\right)-s\left(X_{u}\right)}du,
\end{equation}where $\left(\widetilde{B}_{t}\right)$ is an
$\mathcal{F}_{t}^{\sigma\left(\overline{X}_{T_{0}}\right)}$ Brownian
Motion, stopped at $T_{0}$ and independent of
$\overline{X}_{T_{0}}$. Consequently, we have the following path
decomposition: conditionally on $\overline{X}_{T_{0}}=m$:
\begin{enumerate}
\item the process $\left(X_{t};\;t\leq g\right)$ is a diffusion process started from $x>0$, considered up to $T_{m}$, the first time when it
hits $m$, with infinitesimal generator
$$\frac{1}{2}\dfrac{d^{2}}{dx^{2}}+\left(b\left(x\right)+\dfrac{s'\left(x\right)}{s\left(x\right)}\right)\dfrac{d}{dx}.$$
\item the process $\left(X_{g+t};\;t\leq
T_{0}-g\right)$ is a $\left(\mathcal{F}_{g+t}\right)$ diffusion
process, started from $m$, considered up to $T_{0}$, the first time
when it hits $0$, and is independent of $\left(X_{t};\;t\leq
g\right)$; its infinitesimal generator is given by:
$$\frac{1}{2}\dfrac{d^{2}}{dx^{2}}+\left(b\left(x\right)+\dfrac{s'\left(x\right)}{s\left(x\right)-s\left(m\right)}\right)\dfrac{d}{dx}.$$
\item $\overline{X}_{T_{0}}$ follows the same law as
$s^{-1}\left(\dfrac{1}{U}\right)$, where $U$ follows the uniform law
on $(0,1)$.
\end{enumerate}
\end{itemize}
\end{prop}
\begin{proof}
The proof is exactly the same as the proof of Proposition
\ref{madeco}, so we will not reproduce it here.
\end{proof}

\subsubsection{Geometric Brownian Motion with negative drift} Let $$N_{t}\equiv
\exp\left(2\nu B_{t}-2\nu^{2}t\right),$$where $\left(B_{t}\right)$
is a standard Brownian Motion, and $\nu>0$. With the notation of
Theorem \ref{decoinitial}, we have:
$$S_{t}=\exp\left(\sup_{s\leq t}2\nu\left(
B_{s}-\nu s\right)\right),$$and
$$g=\sup\left\{t:\;\left(
B_{t}-\nu t\right)=\sup_{s\geq 0}\left( B_{s}-\nu
s\right)\right\}.$$ Before stating our proposition, let us mention
that we could have worked with more general continuous exponential
local martingales, but we preferred to keep the discussion as simple
as possible (the proof for more general cases is exactly the same).
\begin{prop}
With the assumptions and notations used above, we have:
\begin{enumerate}
\item The variable $\sup_{s\geq 0}\left(
B_{s}-\nu s\right)$ follows the exponential law of parameter $2\nu$.
\item Every local martingale $X$ is an
$\left(\mathcal{F}_{t}^{\sigma\left(S_{\infty}\right)}\right)$
semimartingale and decomposes as:$$X_{t}=\widetilde{X}_{t}+2\nu
<X,B>_{t\wedge
g}-2\nu\int_{g}^{t}\dfrac{N_{s}}{S_{\infty}-N_{s}}d<X,B>_{s},$$where
$\widetilde{X}_{t}$ is an
$\left(\mathcal{F}_{t}^{\sigma\left(S_{\infty}\right)}\right)$ local
martingale.
\item Conditionally on $S_{\infty}=m$, the process $\left(B_{t}-\nu t;\;t\leq
g\right)$ is a Brownian Motion with  drift $+\nu$ up to the first
hitting time of its maximum $m/2\nu$.
\end{enumerate}
\end{prop}
\begin{proof}
From Doob's maximal equality, $\left(\exp\left(\sup_{s\leq
g}\left(2\nu B_{s}-2\nu^{2}s\right)\right)\right)^{-1}$ follows the
uniform law and hence $\sup_{s\geq 0}\left(B_{s}-\nu s\right)$
follows the exponential law of parameter $2\nu$.

The decomposition formula is a consequence of Theorem
\ref{decoinitial} and the fact that: $dN_{t}=2\nu N_{t}dB_{t}$.

To show $(3)$, it suffices to notice that $B_{t}-\nu t$ is equal to
$\widetilde{B}_{t}+\nu t$ in the filtration
$\left(\mathcal{F}_{t}^{\sigma\left(S_{\infty}\right)}\right)$, with
$\left(\widetilde{B}_{t}\right)$ an
$\left(\mathcal{F}_{t}^{\sigma\left(S_{\infty}\right)}\right)$
Brownian Motion which is independent of $S_{\infty}$.
\end{proof}
\subsubsection{General transient diffusions} Now, we consider
$\left(R_{t}\right)$, a transient diffusion with values in
$\left[0,\infty\right)$, which has $\left\{0\right\}$ as entrance
boundary. Let $s$ be a scale function for $R$, which we can choose
such that: $$s\left(0\right)=-\infty, \text{ and }
s\left(\infty\right)=0.$$ Then, under the law $\mathbf{P}_{x}$, for
any $x>0$, the local martingale
$\left(N_{t}=\dfrac{s\left(R_{t}\right)}{s\left(x\right)},\;t\geq
0\right)$ satisfies the conditions of Theorem \ref{decoinitial}, and
we have:$$\mathbf{P}_{x}\left(g>
t|\mathcal{F}_{t}\right)=\dfrac{s\left(R_{t}\right)}{s\left(I_{t}\right)}\,$$where
$$g=\sup\left\{t:\; R_{t}=I_{t}\right\},$$and
$$I_{t}=\inf_{s\leq t}R_{s}.$$We thus recover results of Jeulin
(\cite{jeulin}, Proposition 6.29, p.112) by other means. Jeulin used
this formula and gave a quick proof of  a theorem of David Williams
(\cite{williams2}), using initial enlargement of filtrations
arguments. Our proof would follow the same lines and so we refer to
the book of Jeulin. We would rather detail an interesting example:
the three dimensional Bessel process.
\begin{prop}
Let $\left(R_{t}\right)$ be a three dimensional Bessel process
starting from $1$, and set, as above, $I_{t}=\inf_{s\leq t}R_{s}$,
and $g=\sup\left\{t:\; R_{t}=I_{t}\right\}$. Define $\rho$ by:
$\rho= \sup\left\{t<g:\;\dfrac{I_{t}}{R_{t}}=\inf_{u\leq
g}\dfrac{I_{u}}{R_{u}}\right\}.$Then:
\begin{enumerate}
\item The variable $\dfrac{I_{\rho}}{R_{\rho}}$ follows the uniform law on $(0,1)$ and, conditionally on $I_{\rho}=rR_{\rho}$, $\left(R_{t},t\leq T_{r}\right)$ is a three
dimensional Bessel process starting from $1$, up to the first time
$T_{r}$ when $I_{t}=rR_{t}$.
\item $I_{\infty}\equiv I_{g}$ follows the uniform law on $(0,1)$;
\item Conditionally on $I_{\infty}=r$, the process $\left(R_{t},\;t\leq
g\right)$ is a Brownian Motion starting from $1$ and stopped when it
first hits $r$.
\end{enumerate}
\end{prop}
\begin{proof}
There exists $\left(\beta\right)_{t\geq 0}$, a Brownian Motion, such
that
$$R_{t}=1+\beta_{t}+\int_{0}^{t}\dfrac{ds}{R_{s}}.$$

$(1)$ follows easily from the results of Subsection
\ref{secpta}.
Now, from Ito's formula, it follows that
$$\dfrac{1}{R_{t}}=1-\int_{0}^{t}\dfrac{d\beta_{s}}{R_{s}^{2}};$$hance, it is a local
martingale. In
$\left(\mathcal{F}_{t}^{\sigma\left(I_{\infty}\right)}\right)$,
$$\beta_{t\wedge g}=\widetilde{\beta}_{t}-\int_{0}^{t\wedge
g}\dfrac{ds}{R_{s}},$$where $\left(\widetilde{\beta}_{t}\right)$ is
an $\left(\mathcal{F}_{t}^{\sigma\left(I_{\infty}\right)}\right)$
Brownian Motion independent of $I_{\infty}$. Hence, $R_{t\wedge g}$
decomposes as
$$R_{t\wedge g}=\widetilde{\beta}_{t}$$ in
$\left(\mathcal{F}_{t}^{\sigma\left(I_{\infty}\right)}\right)$, and
this completes the proof for $(3)$, and $(2)$ is an immediate
consequence of Doob's maximal identity.
\end{proof}
\begin{rem}
The previous method applies to any transient diffusion
$\left(R_{t}\right)_{t\geq 0}$, with values in
$\left(0,\infty\right)$, and which satisfies:
$$R_{t}=x+B_{t}+\int_{0}^{t}duc\left(R_{u}\right),$$where
$c:\mathbb{R}_{+}\rightarrow\mathbb{R}$ allows uniqueness in law for
this equation. These diffusions were studied in
\cite{saichotanemura} to obtain some extension of Pitman's theorem
(see also \cite{zurich}).
\end{rem}
\subsection{Some examples of $Z_{t}$ with jumps}
We shall conclude this paper by giving some explicit examples of
discontinuous $Z's$. Let $X$ be a Poisson process with parameter $c$
and let $N_{t}=X_{t}-ct$. $N$ is a martingale in the natural
filtration $\left(\mathcal{F}_{t}\right)$ of $X$. Every local
martingale $Y$ in this filtration may be written as:
$$Y_{t}=Y_{0}+\int_{0}^{t}k_{s}dN_{s},$$ where $k$ is an
$\left(\mathcal{F}_{t}\right)$ predictable process. Now, for
$f:\mathbb{R}_{+}\rightarrow\mathbb{R}_{+}$ a locally bounded and
Borel function, let
$$\mathcal{E}_{t}^{f}=\exp\left(-\int_{0}^{t}f\left(s\right)dX_{s}+c\int_{0}^{t}\left(1-\exp\left(-f\left(s\right)\right)\right)ds\right)$$
$\mathcal{E}_{t}^{f}$ is an $\mathcal{F}_{t}$ local martingale which
can be represented as:
$$\mathcal{E}_{t}^{f}=1+\int_{0}^{t}\mathcal{E}_{s-}^{f}\left(\exp\left(-f\left(s\right)\right)-1\right)dN_{s}.$$
If $\int_{0}^{\infty}f\left(s\right)ds=\infty$, then
$\lim_{t\rightarrow\infty}\mathcal{E}_{t}^{f}=0$.
\begin{prop}
Let $f$ be a nonnegative locally bounded and Borel function on
$\mathbf{R}_{+}$, such that
$\lim_{t\rightarrow\infty}\mathcal{E}_{t}^{f}=0$.
Define:$$g=\sup\left\{t:\;\mathcal{E}_{t}^{f}=\overline{\mathcal{E}}_{t}^{f}\right\},$$where$$\overline{\mathcal{E}}_{t}^{f}=
\sup_{s\leq t}\mathcal{E}_{s}^{f}.$$ Then:
\begin{enumerate}
\item $\sup_{s\geq
0}\left(-\int_{0}^{t}f\left(s\right)dX_{s}+c\int_{0}^{t}\left(1-\exp\left(-f\left(s\right)\right)\right)ds\right)$
is distributed as a random variable with the exponential law with
parameter $1$;
\item The supermartingale $Z_{t}^{g}$ associated with $g$ is given
by:
$$ \mathbf{P}\left(g>t\mid\mathcal{F}_{t}\right)=\dfrac{\mathcal{E}_{t}^{f}}{\overline{\mathcal{E}}_{t}^{f}};$$
\item Every $\mathcal{F}_{t}$ local martingale $Y_{t}\left(=\int_{0}^{t}k_{s}dN_{s}\right)$ is a
semimartingale in the filtration
$\mathcal{F}_{t}^{\sigma\left(\overline{\mathcal{E}}_{\infty}^{f}\right)}$,
with canonical
decomposition:$$Y_{t}=\widetilde{Y}_{t}+c\int_{0}^{t\wedge
g}k_{s}\left(\exp\left(-f\left(s\right)\right)-1\right)ds-c\int_{g}^{t}k_{s}\left(\exp\left(-f\left(s\right)\right)-1\right)\dfrac{\mathcal{E}_{s}^{f}}{\overline{\mathcal{E}}_{\infty}^{f}-\mathcal{E}_{s}^{f}}ds,$$
where $\widetilde{Y}_{t}$ is an
$\mathcal{F}_{t}^{\sigma\left(\overline{\mathcal{E}}_{\infty}^{f}\right)}$
local martingale.

\end{enumerate}
\end{prop}

\end{document}